\documentclass[11pt,english,reqno]{amsart}
\usepackage[T1]{fontenc}
\usepackage{verbatim}
\usepackage{amssymb}

\makeatletter
 \theoremstyle{plain}
\newtheorem{thm}{Theorem}[section]
  \theoremstyle{plain}
  \newtheorem{lem}[thm]{Lemma}
  \theoremstyle{remark}
  \newtheorem{conclusion}[thm]{Conclusion}
  \theoremstyle{remark}
  \newtheorem{rem}[thm]{Remark}

\usepackage{amscd}
\usepackage{mathrsfs}
\usepackage{hyperref}
 \numberwithin{equation}{section} 
 \numberwithin{figure}{section} 

\usepackage{babel}
\makeatother
\begin{document}

\newcommand{\fall}{,\quad\text{for all}\quad}
\newcommand{\reals}{\mathbb{R}}
\newcommand{\rthree}{\reals^{3}}
\newcommand{\rsix}{\reals^{6}}
\newcommand{\les}{\leqslant}
\newcommand{\ges}{\geqslant}
\newcommand{\dee}{\mathrm{\mathrm{d}}}
\newcommand{\from}{\colon}
\newcommand{\tto}{\longrightarrow}
\newcommand{\abs}[1]{\left|#1\right|}
\newcommand{\isom}{\cong}
\newcommand{\comp}{\circ}
\newcommand{\cl}[1]{\overline{#1}}
\newcommand{\fun}{\varphi}
\newcommand{\interior}{\mathrm{Int\,}}
\newcommand{\sign}{\mathrm{sign\,}}
\newcommand{\dimension}{\mathrm{dim\,}}
\newcommand{\esssup}{\mathrm{ess}\,\sup}
\newcommand{\ess}{\mathrm{{ess}}}
\newcommand{\kernel}{\mathrm{Kernel\,}}
\newcommand{\support}{\mathrm{Supp}\,}
\newcommand{\image}{\mathrm{Image}}
\newcommand{\resto}[1]{|_{#1}}
\newcommand{\incl}{\iota}
\newcommand{\rest}{\rho}
\newcommand{\extnd}{e_{0}}
\newcommand{\proj}{\pi}
\newcommand{\ino}[1]{\int\limits _{#1}}
\newcommand{\half}{\frac{1}{2}}
\newcommand{\shalf}{{\scriptstyle \half}}
\newcommand{\third}{\frac{1}{3}}
\newcommand{\empt}{\varnothing}
\newcommand{\paren}[1]{\left(#1\right)}
\newcommand{\bigp}[1]{\bigl(#1\bigr)}
\newcommand{\braces}[1]{\left\{  #1\right\}  }
\newcommand{\sqbr}[1]{\left[#1\right]}
\newcommand{\norm}[1]{\|#1\|}
\newcommand{\pis}{x}
\newcommand{\pib}{X}
\newcommand{\body}{B}
\newcommand{\bdry}{\partial}
\newcommand{\gO}{\varOmega}
\newcommand{\reg}{\gO}
\newcommand{\bdom}{\bdry\gO}
\newcommand{\bndo}{\partial\gO}
\newcommand{\pbndo}{\Gamma}
\newcommand{\bndoo}{\pbndo_{0}}
 \newcommand{\bndot}{\pbndo_{t}}
\newcommand{\cloo}{\cl{\gO}}
\newcommand{\nor}{\nu}
\newcommand{\dA}{\,\dee A}
\newcommand{\dV}{\,\dee V}
\newcommand{\eps}{\varepsilon}
\newcommand{\vect}{v}
\newcommand{\vs}{\mathbf{W}}
\newcommand{\avs}{\mathbf{V}}
\newcommand{\vbase}{e}
\newcommand{\vf}{w}
\newcommand{\avf}{u}
\newcommand{\stn}{\varepsilon}
\newcommand{\rig}{r}
\newcommand{\rigs}{\mathscr{R}}
\newcommand{\qrigs}{\!/\!\rigs}
\newcommand{\dis}{\chi}
\newcommand{\fc}{F}
\newcommand{\st}{\sigma}
\newcommand{\bfc}{b}
\newcommand{\sfc}{t}
\newcommand{\stm}{S}
\newcommand{\sts}{\varSigma}
\newcommand{\ebdfc}{T}
\newcommand{\optimum}{\st^{\mathrm{opt}}}
\newcommand{\scf}{K}
\newcommand{\cee}[1]{C^{#1}}
\newcommand{\lone}{L^{1}}
\newcommand{\linf}{L^{\infty}}
\newcommand{\ofbdo}{(\bndo)}
\newcommand{\ofclo}{(\cloo)}
\newcommand{\vono}{(\gO,\rthree)}
\newcommand{\vonbdo}{(\bndo,\rthree)}
\newcommand{\vonbdoo}{(\bndoo,\rthree)}
\newcommand{\vonbdot}{(\bndot,\rthree)}
\newcommand{\vonclo}{(\cl{\gO},\rthree)}
\newcommand{\strono}{(\gO,\reals^{6})}
\newcommand{\sob}{W_{1}^{1}}
\newcommand{\sobb}{\sob(\gO,\rthree)}
\newcommand{\lob}{\lone(\gO,\rthree)}
\newcommand{\lib}{\linf(\gO,\reals^{12})}
\newcommand{\ofO}{(\gO)}
\newcommand{\oneo}{{1,\gO}}
\newcommand{\onebdo}{{1,\bndo}}
\newcommand{\info}{{\infty,\gO}}
\newcommand{\infclo}{{\infty,\cloo}}
\newcommand{\infbdo}{{\infty,\bndo}}
\newcommand{\ld}{LD}
\newcommand{\ldo}{\ld\ofO}
\newcommand{\ldoo}{\ldo_{0}}
\newcommand{\trace}{\gamma}
\newcommand{\pr}{\proj_{\rigs}}
\newcommand{\pq}{\proj}
\newcommand{\qr}{\,/\,\reals}
\newcommand{\aro}{S_{1}}
\newcommand{\art}{S_{2}}
\newcommand{\mo}{m_{1}}
\newcommand{\mt}{m_{2}}
\newcommand{\yieldc}{B}
\newcommand{\yieldf}{Y}
\newcommand{\trpr}{\pi_{P}}
\newcommand{\devpr}{\pi_{\devsp}}
\newcommand{\prsp}{P}
\newcommand{\devsp}{D}
\newcommand{\ynorm}[1]{\|#1\|_{\yieldf}}
\newcommand{\colls}{\Psi}
\newcommand{\ssx}{S}
\newcommand{\smap}{s}
\newcommand{\smat}{\chi}
\newcommand{\sx}{e}
\newcommand{\snode}{P}
\newcommand{\node}{p}
\newcommand{\elem}{e}
\newcommand{\nel}{L}
\newcommand{\el}{l}
\newcommand{\ipln}{\phi}
\newcommand{\ndof}{D}
\newcommand{\dof}{d}
\newcommand{\nldof}{N}
\newcommand{\ldof}{n}
\newcommand{\lvf}{\chi}
\newcommand{\lfc}{\varphi}
\newcommand{\amat}{A}
\newcommand{\snomat}{E}
\newcommand{\femat}{E}
\newcommand{\tmat}{T}
\newcommand{\fvec}{f}
\newcommand{\snsp}{\mathcal{S}}
\newcommand{\slnsp}{\Phi}
\newcommand{\ro}{r_{1}}
\newcommand{\rtwo}{r_{2}}
\newcommand{\rth}{r_{3}}

\title{Load Capacity of Bodies}

\author{Reuven Segev}

\curraddr{Department of Mechanical Engineering\\
Ben-Gurion University\\
Beer-Sheva, Israel}

\email{rsegev@bgu.ac.il}

\urladdr{http://www.bgu.ac.il/\textasciitilde{}rsegev}

\keywords{Continuum mechanics, stress analysis, optimization, plasticity, trace.}

\thanks{\today }

\dedicatory{\emph{}}

\subjclass{35Q72; 46E35}

\begin{abstract}
For the stress analysis in a plastic body $\gO$, we prove that there
exists a maximal positive number $C$, the \emph{load capacity ratio,}
such that the body will not collapse under any external traction field
$\sfc$ bounded by $C\yieldf_{0}$, where $\yieldf_{0}$ is the elastic
limit.  The load capacity ratio depends only on the geometry of the
body and is given by\[
\frac{1}{C}=\sup_{\vf\in\ldo_{\devsp}}\frac{\int_{\bndot}\abs{\vf}\dA}{\int_{\gO}\abs{\eps(\vf)}\dV}=\norm{\trace_{\devsp}}.\]
Here, $\ldo_{\devsp}$ is the space of isochoric vector fields $\vf$
for which the corresponding stretchings $\eps(\vf)$ are assumed to
be integrable and $\trace_{\devsp}$ is the trace mapping assigning
the boundary value $\trace_{\devsp}(\vf)$ to any $\vf\in\ldo_{\devsp}$.

\end{abstract}
\maketitle

\section{Introduction}


\noindent Consider a homogeneous isotropic elastic-perfectly plastic
body $\gO$. We prove in this study that there exists a maximal positive
number $C$, to which we refer as the \emph{load capacity ratio,}
such that the body will not collapse under any external traction field
$\sfc$ bounded by $C\yieldf_{0}$, where $\yieldf_{0}$ is the elastic
limit. Thus, while the limit analysis factor of the theory of plasticity
(e.g., \cite{Chris1,Chris2}) pertains to a specific distribution
of external loading, the load capacity ratio is independent of the
distribution of external loading and implies that no collapse will
occur for any field $\sfc$ on $\bndo$ as long as\begin{equation}
\ess\sup_{y\in\bndo}\abs{\sfc(y)}<C\yieldf_{0}.\end{equation}
Collapse will occur for some $\sfc$ if the bound is $C'\yieldf_{0}$
with any $C'>C$.

The load capacity ratio depends only on the geometry of the body and
we prove below that\begin{equation}
\frac{1}{C}=\sup_{\vf\in\ldo_{\devsp}}\frac{\int_{\bndot}\abs{\vf}\dA}{\int_{\gO}\abs{\eps(\vf)}\dV}=\norm{\trace_{\devsp}}.\end{equation}
Here, $\ldo_{\devsp}$ is the space of isochoric integrable vector
fields $\vf$ for which the corresponding stretchings $\eps(\vf)$
are assumed to be integrable and $\trace_{\devsp}$ is the trace mapping
assigning the boundary value $\trace_{\devsp}(\vf)$ to any $\vf\in\ldo_{\devsp}$.

The notion of load capacity ratio is an application to plasticity
 of ideas presented in our previous work \cite{segev03,segev04,PeretzSegev05}
where we consider stress fields on bodies whose maxima are the least.
The general setting may be described as follows.

Let $\gO$ represent the region occupied by the body in space so the
body is supported on a part $\bndoo$ of its boundary and let $\sfc$
be the external surface traction acting on the part $\bndot$ of its
boundary. Body forces may be included in the analysis using the same
methods as in \cite{segev04} but for the sake of simplicity we omit
them here. The mechanical properties of the body are not specified,
and so, there is a class of stress fields that satisfy the equilibrium
conditions with the external loading. (Clearly, distinct distributions
of the mechanical properties within the body will result in general
distinct equilibrating stress distributions.) Each equilibrating stress
field in this class has its own maximal value, and we denote by $\optimum_{\sfc}$
the least maximum.

Specifically, the magnitude of the stress field at a point is evaluated
using a norm on the space of matrices. It is noted that yield conditions
in plasticity usually use semi-norms on the space of stress matrices
rather than norms. By the maximum of a stress field we mean the essential
supremum over the body of its magnitude. Thus, we ignore excessive
values on regions of zero volume. The traction fields that we admit
are essentially bounded also. The set $\gO$ is assumed to be open,
bounded and its boundary is assumed to be smooth. Furthermore, it
is assumed that $\bndot$ and $\bndoo$ are disjoint open subsets
of the boundary whose closures cover the boundary, and that their
closures intersect on a smooth curve.

Subject to these assumptions (see further details in Section \ref{sec:Notation-and-Preliminaries})
our first result is

\begin{thm}
\label{thm:1}\textup{(}i\textup{)} \textbf{The Existence of stresses}.
Given an essentially bounded traction field $\sfc$ on $\bndot$,
there is a collection $\sts_{\sfc}$ of essentially bounded symmetric
tensor fields, interpreted physically as stress fields, that represent
$\sfc$ in the form\begin{equation}
\ino{\bndot}\sfc\cdot\vf\dA=\ino{\gO}\st_{ij}\eps_{ij}(\vf)\dV,\quad\text{for all}\quad\st\in\sts_{\sfc},\ \vf\in\cee{\infty}\vonclo,\label{eq:PrincVirtWork}\end{equation}
where, $\eps(\vf)=\half(\nabla\vf+\nabla\vf^{T})$.

\noindent \textup{(}ii\textup{)} \textbf{The Existence of optimal
stress fields}. There is a stress field $\widehat{\st}\in\sts_{\sfc}$
such that\begin{equation}
\optimum_{\sfc}=\inf_{\st\in\sts_{\sfc}}\braces{\ess\sup_{\pis\in\gO}\abs{\st(\pis)}}=\ess\sup_{\pis\in\gO}\abs{\widehat{\st}(x)}.\end{equation}
\textup{(}iii\textup{)} \textbf{The expression for} \textup{$\optimum_{\sfc}$}.
The optimum satisfies\begin{equation}
\optimum_{\sfc}=\sup_{\vf\in\cee{\infty}\vonclo}\frac{\abs{\int_{\bndot}\sfc\cdot\vf\dA}}{\int_{\gO}\abs{\eps(\vf)}\dV}\,,\label{eq:EpressForOptim}\end{equation}
where the magnitude of $\eps(\vf)(x)$ is evaluated using the norm
dual to the one used for the values of stresses.
\end{thm}
Item (\emph{i}) above is of theoretical interest. It is a representation
theorem for the virtual work performed by the traction field using
tensor fields that we naturally interpret as stresses. It should be
noted that the existence of stress is not assumed here a-priori. The
expression for the representation by stresses turns out to be the
principle of virtual work (\ref{eq:PrincVirtWork}). Thus, the equilibrium
conditions are derived mathematically on the basis of quite general
assumptions.  Item (\emph{i}) also ensures us that the representing
stress fields are also essentially bounded. Item (\emph{ii}) states
that the optimal value is actually attainable for some stress field
and not just as a limit process.

Next, we consider \emph{generalized stress concentration factors}
for the given body. For a given external loading, traditional stress
concentration factors are used by engineers to specify the ratio between
the maximal stress in the body and the maximum nominal stress obtained
using simplified formulas where various geometric irregularities are
not taken into account. Regarding these nominal stresses as boundary
traction fields, we formulate the notion of a stress concentration
factor for a stress field $\st$ in equilibrium with the traction
$\sfc$ mathematically as the ratio between the maximal stress and
the maximum traction. Specifically, we set\begin{equation}
\scf_{\sfc,\st}=\frac{\ess\sup_{\pis\in\gO}\abs{\st(\pis)}}{\ess\sup_{y\in\bndot}\abs{\sfc(y)}}\,.\end{equation}
In particular, the optimal stress concentration factor for the given
traction $\sfc$, is\begin{equation}
\scf_{\sfc}=\inf_{\st\in\sts_{\sfc}}\braces{\scf_{\sfc,\st}}=\frac{\optimum_{\sfc}}{\ess\sup_{y\in\bndot}\abs{\sfc(y)}}\,.\end{equation}
Finally, realizing that engineers may be uncertain as to the nature
of the external loading, we let the external loading vary and define
the \emph{generalized stress concentration factor}, a purely geometric
property of the body $\gO$, as\begin{equation}
\scf=\sup_{\sfc}\braces{\scf_{\sfc}},\end{equation}
where $\sfc$ varies over all essentially bounded traction fields.
In other words, $K$ is the worst possible optimal stress concentration
factor. Using the result on optimal stresses, we prove straightforwardly
the following

\begin{thm}
\label{thm:SCF-2}The generalized stress concentration factor satisfies\begin{equation}
K=\sup_{\vf\in\cee{\infty}\vonclo}\frac{\int_{\bndot}\abs{\vf}\dA}{\int_{\gO}\abs{\eps(\vf)}\dV}=\norm{\trace_{0}},\label{eq:ThmSCF}\end{equation}
where $\trace_{0}$ is the trace mapping for vector fields satisfying
the boundary conditions on $\bndoo$.%
\footnote{Further details on $\trace_{0}$ are described in Section \ref{sec:Constructions-Boun-Cond}.%
}
\end{thm}
To prove the theorems we use standard tools of analysis and the theory
of $\ld$-spaces given by \cite{Temam81,Temam85,TemamStrang1}. Results
analogous to Theorems (\ref{thm:1}) and (\ref{thm:SCF-2}) were presented
in our earlier work cited above. In \cite{segev03}, a weaker form
of equilibrium is assumed, and in all earlier work we did not consider
boundary conditions for the displacements on $\bndoo$.

Next we turn to the adaptation needed for the application to plasticity.
It is assumed that the yield function is a norm on the space of matrices
applied to the deviatoric component of the stress matrix. Thus, it
turns out that the same mathematical structure applies if we consider
isochoric (incompressible) vector fields in the suprema of Equations
(\ref{eq:EpressForOptim}) and (\ref{eq:ThmSCF}). For example, the
analog of Equation (\ref{eq:EpressForOptim}) is\begin{equation}
\optimum_{\sfc}=\sup_{\vf\in\ldo_{\devsp}}\frac{\abs{\int_{\bndot}\sfc\cdot\vf\dA}}{\int_{\gO}\abs{\eps(\vf)}\dV}\,,\label{eq:OptStPlast}\end{equation}
where $\ldo_{\devsp}$ is the collection of isochoric vector fields
having integrable strains.

It turns out that optimal stresses are related to limit analysis of
plasticity. In fact, the limit analysis factor $\lambda^{*}$ (see
Remark \ref{rem:RelToChris} and \cite{Chris1,Chris2,TemamStrang2})
is simply given by\begin{equation}
\lambda^{*}=\frac{\yieldf_{0}}{\optimum_{\sfc}}\,.\end{equation}
Furthermore, the expression for the optimal stress of Equation (\ref{eq:OptStPlast})
is implied mathematically by the results of Christiansen and Temam
\& Strang (\cite{Chris1,Chris2,TemamStrang2}). This implies that
the optimal stress fields do not require a special designs of non-homogeneous
material properties but occur for the frequently used models of elastic-plastic
bodies. In particular, elastic-plastic material will attain the optimal
stress field independently of the distribution of the external load.

We take advantage of these observations and introduce here the notion
of load capacity ratio---a purely geometric property of the body.
As described above, the load capacity ratio may be conceived as a
universal limit design factor, which is independent of the distribution
of the external loading. It immediately follows from its definition
that \begin{equation}
C=\frac{1}{K}\,.\end{equation}

Section \ref{sec:Notation-and-Preliminaries} introduces the notation,
assumptions and some background material. In particular, the space
$\ldo$ of vector fields of integrable stretchings (or linear strains)
(see \cite{Temam81,Temam85,TemamStrang1}) is described. Following
some preliminary material concerning the boundary conditions in Section
\ref{sec:Constructions-Boun-Cond}, the proof of the theorems is given
in Section \ref{sec:Proof-of-Thm1}, with some additional details
in Appendix \ref{sec:The-Technical-Lemma}. The adaptation to plasticity
theory, including the introduction of the load capacity ratio, is
presented in Section \ref{sec:Load-Capac-for-Plas}.

It is noted that for structures, i.e., the finite dimensional approximations
of the present setting, the expression for $C$ may be set as a linear
programming problem. Thus, the load capacity ratio and the generalized
stress concentration factor may be computed approximately using standard
algorithms.

I would like to thank anonymous reviewers for pointing out the relation
between the notion of optimal stress and the expression I obtained
for it, and the limit analysis factor and the expressions for it in
the work of Christiansen and Temam \& Strang (\cite{Chris1,Chris2,TemamStrang2}).

\section{Notation and Preliminaries\label{sec:Notation-and-Preliminaries}}

\subsection{Basic variables}

We consider a body under a given configuration in space. The space
is modelled simply by $\rthree$ and the image of the body under the
given configuration is the subset $\gO\subset\rthree$. It is assumed
that $\gO$ is open and bounded and that it has a $\cee{^{1}}$-boundary
$\bndo$. Furthermore, there are two open subsets $\bndoo\subset\bndo$,
and $\bndot\subset\bndo$ such that $\bndoo$ is the region where
the body is supported and $\bndot$ is the region where the body is
not supported so that a surface traction field $\sfc$ may be exerted
on the body on $\bndot$. Thus, it is natural to assume that $\bndoo$
and $\bndot$ are nonempty and disjoint, $\cl{\pbndo}_{0}\cup\cl{\pbndo}_{\sfc}=\bndo$,
and $\Lambda=\bdry\bndoo=\bdry\bndot$ is a differentiable 1-dimensional
submanifold of $\bndo$. (The regularity assumptions may be generalized
without affecting the validity of the constructions below.)

Basic objects in the construction are spaces of generalized velocity
fields. A generic \emph{generalized velocity field} (alternatively,
\emph{virtual velocity} or \emph{virtual displacement}) will be denoted
by $\vf$. In the sequel we consider a number of Banach spaces containing
generalized velocities and a generic space of generalized velocities
will be a denoted by $\vs$. Generalized forces will be elements of
the dual space $\vs^{*}$. Thus, a \emph{generalized force} $\fc$
is a bounded linear functional $\fc\from\vs\tto\reals,$ such that
$\fc(\vf)$ is interpreted as the virtual power (virtual work) performed
by the force for the generalized velocity $\vf$. We recall that the
dual norm of a linear functional $\fc$ is defined as\begin{equation}
\norm{\fc}=\sup_{\vf\in\vs}\frac{\abs{\fc(\vf)}}{\norm{\vf}}\,.\end{equation}

\subsection{Virtual stretchings (linear strains) and stresses}

As an example for the preceding paragraph, consider the case where
$\vs$ is the space $\lone\strono$ of $\lone$-symmetric tensor fields
on $\gO$. A typical element $\eps\in\lone\strono$ is interpreted
as a \emph{virtual stretching field} or a \emph{linear strain field}.
We will use $\abs{\eps(\pis)}$ to denote the norm of the matrix $\eps(\pis)$.
Various such norms are described in \cite{PeretzSegev05}. Thus,\begin{equation}
\norm{\eps}_{1}=\ino{\gO}\abs{\eps(x)}\dV.\end{equation}
The dual space $\lone\strono^{*}=\linf\strono$ contains symmetric
essentially bounded tensor fields $\st$ that act on the stretching
fields by \begin{equation}
\st(\eps)=\ino{\gO}\st(\pis)(\eps(\pis))\dV.\end{equation}
Here, we use the same notation for the functional $\st$ and the essentially
bounded tensor field representing it and we regard the matrix $\st(\pis)$
as a linear form on the space of matrices so $\st(\pis)(\eps(\pis))=\st(\pis)_{ij}\eps(\pis)_{ji}$.
Naturally, an element $\st\in\linf\strono$ is interpreted as a \emph{stress
field}. The dual norm of a stress field is given as\begin{equation}
\norm{\st}=\norm{\st}_{\infty}=\ess\sup_{\pis\in\gO}\abs{\st(\pis)}.\end{equation}
 Here, $\abs{\st(\pis)}$ is calculated using the norm on the space
of matrices which is dual to the one used for the evaluation of $\abs{\eps(\pis)}$
(see \cite{PeretzSegev05} for details). Thus, the choice of the space
$\lone\strono$ for stretchings is natural when one is looking for
the maximum of the stress tensor.

\subsection{The space of boundary velocity fields and boundary tractions}

As another example to be used later, consider the space $\lone\vonbdot$
of integrable vector fields on the {}``free'' part of the boundary.
Its dual space is \begin{equation}
\lone\vonbdot^{*}=\linf\vonbdot,\end{equation}
so that a generalized force in this case will be represented by an
essentially bounded vector field $\sfc$ on $\bndot$. Using the same
notation for the functional and the vector field representing it,
we have\begin{equation}
\sfc(\avf)=\ino{\bndot}\sfc(y)\cdot\avf(y)\dA\end{equation}
so $\sfc$ may be interpreted as a traction field on $\bndot$ as
expected. The dual norm of the traction field $\sfc$ is \begin{equation}
\norm{\sfc}=\norm{\sfc}_{\infty}=\ess\sup_{y\in\bndot}\abs{\sfc(y)},\end{equation}
again, the relevant maximum.

\subsection{The space $\ldo$ and its elementary properties\label{sub:The-space-LDO}}

A central role in the subsequent analysis is played by the space $\ldo$
containing vector fields of integrable stretchings (see \cite{TemamStrang1,Temam81,Temam85}).
We summarize below its definition and basic relevant properties (see
\cite{Temam85} for proofs and details).

\subsubsection{Definition}

For an integrable vector field $\vf\in\lone\vono$, let $\nabla\vf$
denote its distributional gradient and consider the corresponding
stretching (a tensor distribution)\begin{equation}
\eps(\vf)=\half(\nabla\vf+\nabla\vf^{T}).\end{equation}
The vector field $\vf$ has an integrable stretching if the distribution
$\eps(\vf)$ is an integrable symmetric tensor field, i.e., it belongs
to $\lone\strono$. For the sake of simplifying the notation, we use
$\eps$ for both the stretching mapping here and its value in the
example above. The space $\ldo$ is defined by\begin{equation}
\ldo=\braces{\vf\from\gO\to\rthree;\;\vf\in\lone\vono,\;\eps(\vf)\in\lone\strono}.\end{equation}
A natural norm is provided by\begin{equation}
\norm{\vf}=\norm{\vf}_{\ld}=\norm{\vf}_{1}+\norm{\eps(\vf)}_{1}\label{eq:ld-norm}\end{equation}
and it induces on $\ldo$ a Banach space structure. Clearly, the stretching
mapping\begin{equation}
\eps\from\ldo\tto\lone\strono\end{equation}
is linear and continuous.

\subsubsection{Approximations}

The space of restrictions to $\gO$ of smooth mappings in $\cee{\infty}\vonclo$,
is dense in $\ldo$, so any $\ld$-vector field may be approximated
by restrictions of smooth vector fields defined on the closure $\cl{\gO}$.

\subsubsection{Trace mapping}

There is a unique continuous and linear trace mapping\begin{equation}
\trace\from\ldo\tto\lone\vonbdo\end{equation}
satisfying the consistency condition\begin{equation}
\trace(\avf\resto{\gO})=\avf\resto{\bndo}\end{equation}
for any continuous mapping $\avf\in\cee{0}\vonclo$. Furthermore,
the trace mapping is surjective. Thus, although $\ld$-mappings are
defined on the open set $\gO$, they have meaningful $\lone$ boundary
values.

\subsubsection{Equivalent norm\label{sub:Equivalent-norm}}

Let $\pbndo$ be an open subset of $\bndo$ and for $\vf\in\ldo$
let\begin{equation}
\norm{\vf}_{\pbndo}=\ino{\pbndo}\abs{\trace(\vf)}\dA+\norm{\eps(\vf)}_{1},\label{eq:EquivNorm}\end{equation}
then, $\norm{\vf}_{\pbndo}$ is a norm on $\ldo$ which is equivalent
to the original norm defined in Equation (\ref{eq:ld-norm}).

\section{Constructions Associated with the Boundary Conditions\label{sec:Constructions-Boun-Cond}}

\subsection{The space $\lone\vonbdo_{0}$}

Let $\lone\vonbdo_{0}\subset\lone\vonbdo$ be the vector space of
vector fields on $\bndo$ such that for each $\avf\in\lone\vonbdo_{0}$,
$\avf(y)=0$ for almost all $y\in\bndoo$. It is noted that the restriction
mapping\begin{equation}
\rest_{0}\from\lone\vonbdo\tto\lone\vonbdoo,\quad\rest_{0}(\avf)=\avf\resto{\bndoo}\end{equation}
is linear and continuous. Thus, since \begin{equation}
\lone\vonbdo_{0}=\rest_{0}^{-1}\{0\},\end{equation}
$\lone\vonbdo_{0}$ is a closed subspace of $\lone\vonbdo$.

The restriction mapping \begin{equation}
\rest_{\sfc}\from\lone\vonbdo_{0}\tto\lone\vonbdot,\quad\rest_{\sfc}(\avf)=\avf\resto{\bndot}\end{equation}
is also linear and continuous. In addition, as $\bdry\bndoo=\bdry\bndot=\Lambda$
have zero area measure,\begin{equation}
\ino{\bndo}\abs{\avf}\dA=\ino{\bndot}\abs{\rest_{\sfc}(\avf)}\dA,\quad\avf\in\lone\vonbdo_{0},\end{equation}
so $\rest_{t}$ is a norm-preserving injection.

Consider the zero extension mapping\begin{equation}
\extnd\from\lone\vonbdot\tto\lone\vonbdo_{0},\ \end{equation}
defined by\begin{equation}
\extnd(\avf)(y)=\begin{cases}
\avf(y) & \text{for\ }y\in\bndot,\\
0 & \text{for\ }y\notin\bndot.\end{cases}\end{equation}
Clearly, $\rest_{\sfc}\comp\extnd$ is the identity on the space $\lone\vonbdot$.
Moreover, for any $\avf\in\lone\vonbdo_{0}$, $\extnd(\rest_{\sfc}(\avf))(y)=\avf(y)$
almost everywhere (except for $y\in\Lambda$), so $\extnd\comp\rest_{\sfc}$
is the identity on $\lone\vonbdo_{0}$. We conclude,

\begin{lem}
\label{lem:extension and restiction}The mappings $\rest_{\sfc}$
and $\extnd$ induce an isometric isomorphism of the spaces $\lone\vonbdo_{0}$
and $\lone\vonbdot$. The dual mappings $\extnd^{*}$ and $\rest_{\sfc}^{*}$
induce an isometric isomorphism of the spaces $\lone\vonbdot^{*}$
and $\lone\vonbdo_{0}^{*}$. Every element $\sfc_{0}\in\lone\vonbdo_{0}^{*}$
is represented uniquely by an essentially bounded $\sfc\in\linf\vonbdot$
in the form\begin{equation}
\sfc_{0}(\avf)=\ino{\bndot}\sfc\cdot\avf\dA.\end{equation}

\end{lem}

\subsection{The space $\ldoo$ of velocity fields satisfying the boundary conditions}

Recalling the definition of the equivalent norm on $\ldo$ in Equation~(\ref{eq:EquivNorm}),
we set $\pbndo=\bndoo$ in that equation. Henceforth, we will use
on $\ldo$ the equivalent norm\begin{equation}
\norm{\vf}=\norm{\vf}_{\bndoo}=\ino{\bndoo}\abs{\trace(\vf)}\dA+\norm{\eps(\vf)}_{1}.\label{eq:UsedNormLdo}\end{equation}

Consider the vector subspace $\ldo_{0}$ defined by\begin{equation}
\ldoo=\trace^{-1}\braces{\lone\vonbdo_{0}}\subset\ldo.\end{equation}
 Thus, $\ldoo$ is the subspace containing vector fields on $\gO$
whose boundary values vanish on $\bndoo$ almost everywhere. Since
$\trace$ is continuous and $\lone\vonbdo_{0}$ is a closed subspace
of $\lone\vonbdo$, $\ldo_{0}$ is a closed subspace of $\ldo$. Combining
this with Lemma (\ref{lem:extension and restiction}) we obtain immediately

\begin{lem}
\label{lem:trace_0}The mapping\begin{equation}
\trace_{0}=\rest_{\sfc}\comp\trace\resto{\ldoo}\from\ldoo\tto\lone\vonbdot\end{equation}
is a linear and continuous surjection. Dually, \begin{equation}
\trace_{0}^{*}=\paren{\trace\resto{\ldoo}}^{*}\comp\rest_{\sfc}^{*}\from\linf\vonbdot\tto\ldoo^{*}\end{equation}
is a continuous injection.
\end{lem}
Observing Equation (\ref{eq:UsedNormLdo}), for each $\vf\in\ldo_{0}$,
\begin{equation}
\norm{\vf}=\norm{\eps(\vf)}_{1}.\label{eq:normForLDOO}\end{equation}

\begin{lem}
\label{lem:EpsInjectiveIsometry}The mapping\begin{equation}
\eps_{0}=\eps\resto{\ldoo}\from\ldoo\to\lone\strono\end{equation}
is an isometric injection.
\end{lem}
\begin{proof}
Equation (\ref{eq:normForLDOO}) implies immediately that $\norm{\vf}=\norm{\eps(\vf)}_{1}$
for all $\vf\in\ldoo$. Being a linear isometry, the zero element
is the only element that is mapped to zero, so $\eps_{0}$ is injective.
In addition to relying on the technical property (\ref{sub:Equivalent-norm})
of $\ldo$ to show that $\eps_{0}$ is injective, it should be mentioned
that this follows from the fact that for any vector field $\vf$ on
$\gO$, $\eps(\vf)=0$ only if $\vf$ is a rigid vector field, i.e.,
if $\vf$ is of the form $\vf(\pis)=a+b\times\pis$, $a,b\in\rthree$.
Now, the only rigid vector field that vanishes on the open set $\bndoo$
is the zero vector field.
\end{proof}

\section{The Mathematical Constructions\label{sec:Proof-of-Thm1}}

Let $\sfc\in\linf\vonbdot$ be a traction field on the free part of
the boundary. Then, $\trace_{0}^{*}(\sfc)$ is an element of $\ldoo^{*}$
representing $\sfc$. The basic properties of elements of $\ldoo^{*}$
are as follows.

\begin{lem}
\label{lem:Optimal Stress}Each $\stm\in\ldoo^{*}$ may be represented
by some non-unique tensor field $\st\in\linf\strono$ in the form\begin{equation}
\stm=\eps_{0}^{*}(\st)\quad\text{or}\quad\stm(\vf)=\ino{\gO}\st(\pis)(\eps_{0}(\vf)(\pis))\dV.\end{equation}
 The dual norm of $\stm$ satisfies\begin{equation}
\norm{\stm}=\inf_{\st}\norm{\st}_{\infty}=\inf_{\st}\braces{\ess\sup_{\pis\in\gO}\abs{\st(\pis)}},\end{equation}
where the infimum is taken is taken over all tensor fields $\st$,
satisfying $\stm=\eps_{0}^{*}(\st)$, i.e, tensor fields representing
$\stm$. There is a $\widehat{\st}\in\linf\strono$ for which the
infimum is attained.
\end{lem}
\begin{proof}
The assertion follows from the fact that $\eps_{0}$ is a linear and
isometric injection as in Lemma (\ref{lem:EpsInjectiveIsometry})
and using the Hahn-Banach theorem. See Appendix~\ref{sec:The-Technical-Lemma}
for the details of the technical lemma used and its proof.
\end{proof}
Applying this lemma to $\stm=\trace_{0}^{*}(\sfc)$ one may draw the
following conclusions.

\begin{conclusion}
Forces on the body given by essentially bounded surface tractions
are represented by tensor fields on the body. These tensor fields
are naturally interpreted as stress fields. The condition that a stress
tensor field $\st$ represents the surface traction $\sfc$ is \begin{equation}
\trace_{0}^{*}(\sfc)=\eps_{0}^{*}(\st),\end{equation}
and explicitly,\begin{equation}
\ino{\bndot}\sfc\cdot\trace_{0}(\vf)\dA=\ino{\gO}\st(\eps_{0}(\vf))\dV,\end{equation}
for each vector field $\vf\in\ldoo$, i.e., a vector field of integrable
stretching satisfying the boundary condition on $\bndoo$. This condition
is just the principle of virtual work which is a weak form of the
equation of equilibrium and the corresponding boundary conditions.
Thus, we have derived both the existence of stresses and the equilibrium
conditions analytically under mild assumptions.

It is noted that the subscript 0, only indicating the restriction
of the various operations to fields satisfying the boundary conditions,
may be omitted above. Also, as the restrictions of smooth vector fields
on $\cl{\gO}$ are dense in $\ldo$, it is sufficient to verify that
the condition holds for smooth fields on $\cl{\gO}$. For such fields,
the integrand on the left may be replaced simply by $\sfc\cdot\vf$.
\end{conclusion}
{}

\begin{conclusion}
There is an optimal stress field $\widehat{\st}$ representing $\sfc$
and\begin{equation}
\norm{\trace_{0}^{*}(\sfc)}=\norm{\widehat{\st}}_{\infty}=\inf_{\st}\braces{\ess\sup_{\pis\in\gO}\abs{\st(\pis)}},\end{equation}
where the infimum is taken over all stress fields $\st$ satisfying
$\trace_{0}^{*}(\sfc)=\eps_{0}^{*}(\st)$, i.e., all stress fields
in equilibrium with $\sfc$. Thus, the infimum on the right is the
optimal maximal stress $\optimum_{\sfc}$. In addition, by the definition
of the dual norm we have\begin{align}
\norm{\trace_{0}^{*}(\sfc)} & =\sup_{\vf\in\ldoo}\frac{\abs{\trace_{0}^{*}(t)(\vf)}}{\norm{\vf}}\\
 & =\sup_{\vf\in\ldoo}\frac{\abs{t(\trace_{0}(\vf))}}{\norm{\eps(\vf)}_{1}}\,,\end{align}
where in the last line we used Equation (\ref{eq:normForLDOO}). We
conclude that\begin{equation}
\optimum_{\sfc}=\sup_{\vf\in\ldoo}\frac{\abs{\int_{\bndot}t\cdot\trace_{0}(\vf)\dA}}{\int_{\gO}\abs{\eps(\vf)}\dV}\,.\label{eq:OptimumResult}\end{equation}
Recalling that the restrictions of smooth mappings on $\cl{\gO}$
are dense in $\ldo$ and that for such mappings the trace mapping
is just the restriction, the optimal stress may be evaluated as\begin{equation}
\optimum_{\sfc}=\sup_{\vf}\frac{\abs{\int_{\bndot}t\cdot\vf\dA}}{\int_{\gO}\abs{\eps(\vf)}\dV}\,,\label{eq:OptimumResult-C-Inf}\end{equation}
where the supremum is taken over all smooth mappings in $\cee{\infty}\vonclo$
that vanish on $\bndoo$.

It is noted that the value of $\optimum_{\sfc}$ depends on the norm
used for strain matrices.
\end{conclusion}
We now turn to the simple proof of Theorem (\ref{thm:SCF-2}).

\begin{proof}
We had\[
\optimum_{\sfc}=\sup_{\vf\in\ldoo}\frac{\abs{\sfc(\trace_{0}(\vf))}}{\norm{\eps(\vf)}_{1}}=\norm{\trace_{0}^{*}(\sfc)}\,,\]
so

\begin{align}
\scf & =\sup_{\sfc\in\linf\vonbdot}\frac{\optimum_{\sfc}}{\norm{\sfc}}=\sup_{\sfc\in\linf\vonbdot}\braces{\frac{\norm{\trace_{0}^{*}(\sfc)}}{\norm{\sfc}}}=\norm{\trace_{0}^{*}}=\norm{\trace_{0}}\end{align}
where the last equality is the standard equality between the norm
of a mapping and the norm of its dual (e.g., \cite[pp., 191-192]{Taylor}).
\end{proof}

\section{Load Capacity for Plastic Bodies\label{sec:Load-Capac-for-Plas}}

The analysis we presented in the previous sections may be applied
to the limit analysis of plastic bodies. While in the preceding analysis
the magnitude of the stress at a point was represented by the norm
of the stress matrix, for the analysis of plasticity, the relevant
quantity is the value of the yield function. The yield function is
usually taken as a seminorm on the space of matrices---a norm on the
deviatoric component of the stress. The necessary adaptation is as
follows.

\subsection{Notation and Preliminaries}

We denote by $\trpr$ the projection on the subspace of spherical
matrices $\prsp$, i.e., \begin{equation}
\trpr(m)=\third m_{ii}I,\end{equation}
and by $\devpr$, the projection on the subspace of deviatoric matrices
$\devsp$ so\begin{equation}
\devpr(m)=m_{D}=m-\trpr(m).\end{equation}
Thus, the pair $(\devpr,\trpr)$ makes an isomorphism of the space
of matrices with $\devsp\oplus\prsp$. We will therefore make the
identification $\rsix=\devsp\oplus\prsp$ and $\reals^{6*}=\devsp^{*}\oplus\prsp^{*}$.
We will use the same notation $\abs{\cdot}$ for both the norm on
$\rsix$, whose elements are interpreted as strain values, and the
dual norm on $\reals^{6*}$, whose elements are interpreted as stress
values (although the norms may be different in general). Thus, we
assume that the yield function is of the form \begin{equation}
\yieldf(m)=\abs{\devpr(m)}.\end{equation}
For example, if we take $\abs{\cdot}$ to be the 2-norm on $\reals^{6*}$
we get the Von-Misses yield criterion. In practical terms this means
that the material yields at a point $x$ when $\yieldf(\st(x))=\yieldf_{0}$
for some limiting yield stress value $\yieldf_{0}\in\reals^{+}$.

Thus, we will extend the foregoing discussion to the case where $\yieldf$,
evidently a seminorm, replaces the norm on the space of stress matrices.
For the space of stress fields we will therefore have the seminorm
$\ynorm{\cdot}$ defined by\begin{equation}
\ynorm{\st}=\norm{\devpr\comp\st}_{\infty}=\ess\sup_{x\in\gO}\yieldf(\st(x)).\end{equation}
The expression defining the optimal stress becomes\begin{equation}
\optimum_{\sfc}=\inf_{\trace^{*}(\sfc)=\eps^{*}(\st)}\ynorm{\st}=\inf_{\trace^{*}(\sfc)=\eps^{*}(\st)}\braces{\norm{\devpr\comp\st}_{\infty}}.\end{equation}
The condition for collapse is $\optimum_{\sfc}\ges\yieldf_{0}$ and
we use $\colls$ to denote the collapse manifold, i.e.,\begin{equation}
\colls=\braces{\sfc\mid\optimum_{\sfc}=\yieldf_{0}}.\end{equation}

\begin{rem}
The expression\[
\optimum_{\sfc}=\inf_{\eps^{*}(\st)=\trace^{*}(\sfc)}\ynorm{\st}\]
for the optimal stress may be refomultated as follows. Recalling that
$\yieldf_{0}$ denotes the yield stress, we write $\st=\st_{1}/\lambda$,
$\ynorm{\st_{1}}=\yieldf_{0}$ and noting that $\norm{\st/\norm{\st}_{\yieldf}}_{\yieldf}=1$,
we are looking for \begin{equation}
\optimum_{\sfc}=\inf_{\substack{\eps^{*}(\st_{1}/\lambda)=\trace^{*}(\sfc),\\
\lambda\in\reals^{+},\,\st_{1}\in\bdry B}
}\norm{\st_{1}/\lambda}_{\yieldf}\end{equation}
where $B$ is the ball in $\linf(\gO,\devsp)$ of radius $\yieldf_{0}$
. Thus,\begin{align}
\optimum_{\sfc} & =\inf_{\substack{\eps^{*}(\st_{1}/\lambda)=\trace^{*}(\sfc),\\
\lambda\in\reals^{+},\,\st_{1}\in\bdry B}
}\frac{\yieldf_{0}}{\lambda},\\
\frac{\yieldf_{0}}{\optimum_{\sfc}} & =\sup\braces{\lambda\mid\,\exists\st_{1}\in\bdry B,\,\eps^{*}(\st_{1})=\trace^{*}(\lambda\sfc)}.\end{align}

Clearly, in the last equation $\partial B$ may be replaced by $\cl{B}$
because if we consider $\st$ with $\norm{\st}_{\yieldf}<1$, then,
$\st_{1}=\st/\norm{\st}_{\yieldf}$ is in $\bdry B$ and the corresponding
$\lambda$ will be multiplied by $\norm{\st}_{\yieldf}<1$.

In the last equation the unit ball $B$ contains the (elastic) states
of the material within the yield surface and we are looking for the
largest multiplication of the force so the resulting stress is within
the yield surface. Thus, we are looking for \begin{equation}
\frac{\yieldf_{0}}{\optimum_{\sfc}}=\lambda^{*}=\sup\lambda,\,\exists\st\in B,\,\eps^{*}(\st)=\trace^{*}(\lambda\sfc).\end{equation}
which is the limit analysis factor (e.g., Christiansen \cite{Chris1,Chris2}
and Teman \& Strang \cite{TemamStrang2}).
\end{rem}
The expression defining the generalized stress concentration factor
assumes the form\begin{equation}
K=\sup_{\sfc}\frac{\optimum_{\sfc}}{\norm{\sfc}_{\infty}}=\sup_{\sfc}\inf_{\trace^{*}(\sfc)=\eps^{*}(\st)}\frac{\ynorm{\st}}{\norm{\sfc}_{\infty}}=\sup_{\sfc}\inf_{\trace^{*}(\sfc)=\eps^{*}(\st)}\frac{\norm{\devpr\comp\st}_{\infty}}{\norm{\sfc}_{\infty}}\,.\end{equation}

For the application to plasticity, we use the term \emph{load capacity}
for \emph{$C=1/K$.} Hence,\begin{equation}
C=\frac{1}{\sup_{t}(\optimum_{\sfc}/\norm{\sfc}_{\infty})}=\inf_{\sfc}\frac{\norm{\sfc}_{\infty}}{\optimum_{\sfc}}\,.\end{equation}
For every loading $\sfc$ we set\begin{equation}
\sfc_{\colls}=\frac{\sfc}{\optimum_{\sfc}/\yieldf_{0}}\,,\end{equation}
 so using $\optimum_{\lambda\sfc}=\norm{\trace^{*}(\lambda\sfc)}=\lambda\optimum_{\sfc}$
for any $\lambda>0$, one has\begin{equation}
\optimum_{\sfc_{\colls}}=\yieldf_{0},\quad\frac{\norm{\sfc}_{\infty}}{\optimum_{\sfc}}=\frac{\norm{\sfc_{\colls}\optimum_{\sfc}/\yieldf_{0}}_{\infty}}{\optimum_{\sfc}}=\norm{\sfc_{\colls}}_{\infty}/\yieldf_{0}.\end{equation}
It follows that for any $\sfc$, $\sfc_{\colls}$ belongs to the collapse
manifold $\colls$ and the operation above is a projection onto the
collapse manifold. Thus,\begin{equation}
C=\inf_{\sfc}\frac{\norm{\sfc}_{\infty}}{\optimum_{\sfc}}=\inf_{\sfc_{\colls}\in\colls}\norm{\sfc_{\colls}}_{\infty}/\yieldf_{0},\end{equation}
and indeed, $C\yieldf_{0}=\inf_{\sfc_{\colls}\in\colls}\norm{\sfc_{\colls}}_{\infty}$
is the largest radius of a ball containing only surface forces for
which collapse does not occur.

\subsection{Constructions associated with the extension to plasticity}

The top row of the following commutative diagram describes the various
kinematic mappings we used in the case of a norm on the space of stresses
as considered in the previous sections. The subspace $\lone(\gO,\devsp)$
of $\lone\strono$ contains deviatoric (isochoric) strain fields and
there is a natural projection $\devpr^{\circ}\from\lone\strono\to\lone(\gO,\devsp)$
given by $\devpr^{\circ}(\dis)=\devpr\comp\dis$. The inclusion of
a subspace in a vector space will be generally denoted as $\incl$.
We will also use the notation\begin{equation}
\ldo_{\devsp}=\eps_{0}^{-1}\braces{\lone(\gO,\devsp)}\end{equation}
for the subspace of isochoric $\ld$-vector fields. Thus, using $\eps_{\devsp}$
and $\trace_{\devsp}$ for the restrictions $\eps_{0}\resto{\ldo_{\devsp}}$
and $\trace_{0}\resto{\ldo_{\devsp}}$, respectively, we have the
following commutative diagram.

\begin{equation}
\begin{CD}
\lone\vonbdot  @<{\trace_0}<<      \ld(\gO)_0
                   @>\eps_0>> \lone(\gO,\rsix)
                               \\
@| {}         @AA{\incl}A       {}   @A{\incl}AA\hspace*{-8.5mm}@VV{\devpr^\circ}V\\
                 \lone\vonbdot  @<\trace_\devsp<<  \ldo_\devsp
                   @>\eps_\devsp>>
                             \lone(\gO,\devsp).
\end{CD}
\end{equation}

The dual diagram is

\begin{equation}
\begin{CD}
\linf\vonbdot  @>{\trace_0}^*>>      \ld(\gO)^*_0
                   @<\eps^*_0<< \linf(\gO,\rsix)
                               \\
@| {}         @VV{\incl^*}V       {}   @V{\incl^*}VV\hspace*{-9mm}@AA{\devpr^{\circ*}}A\\
\linf\vonbdot  @>\trace^*_\devsp>>  \ldo_\devsp^*
                   @<\eps^*_\devsp<<
                             \linf(\gO,\devsp).
\end{CD}
\end{equation}

Since $\eps_{\devsp}$ is just a restriction of $\eps_{0}$, it is
still a linear, norm-preserving injection. Thus, the assertion of
Lemma (\ref{lem:Optimal Stress}) and the subsequent conclusions hold
where $\ldo_{\devsp}$, $\eps_{\devsp}$, and $\trace_{\devsp}$ replace
$\ldo_{0}$, $\eps_{0}$, and $\trace_{0}$, respectively. The expression
for the optimal stress for the plasticity analysis is therefore,\begin{equation}
\optimum_{\sfc}=\inf_{\trace^{*}(\sfc)=\eps^{*}(\st)}\ynorm{\st}=\sup_{\vf\in\ldo_{\devsp}}\frac{\abs{\int_{\bndot}\sfc\cdot\vf\dA}}{\int_{\gO}\abs{\eps(\vf)}\dV}\label{eq:OptStrPlastProof}\end{equation}

Finally, the load capacity ratio is given by \begin{equation}
\frac{1}{C}=K=\sup_{t\in\linf\vonbdot}\frac{\optimum_{\sfc}}{\norm{\sfc}_{\infty}}=\sup_{\vf\in\ldo_{\devsp}}\frac{\int_{\bndot}\abs{\vf}\dA}{\int_{\gO}\abs{\eps(\vf)}\dV}=\norm{\trace_{\devsp}}.\end{equation}

\begin{rem}
\label{rem:RelToChris}Our result (\ref{eq:OptStrPlastProof}) for
the optimal stress associated with $\stm=\trace_{\devsp}^{*}(\sfc)$
is\begin{equation}
\optimum_{\stm}=\sup_{\vf\in\ldo_{\devsp}}\frac{\abs{\stm(\vf)}}{\norm{\eps(\vf)}_{\yieldf}}\,.\label{eq:OptS}\end{equation}
For limit analysis in plasticity it is shown by Christiansen \cite{Chris1,Chris2}
and Teman \& Strang \cite{TemamStrang2} that the kinematic version
of the limit load is equivalent to the statical version above, specifically,
\begin{equation}
\lambda^{*}=\inf_{\stm(\vf)=1}\braces{\sup_{\st\in B}\braces{\eps^{*}(\st)(\vf)}}.\end{equation}
We will show that the two expressions are equivalent for the setting
of stress optimization.

Equation (\ref{eq:OptS}) may be rewritten as \begin{align}
\optimum_{\stm} & =\sup_{\vf\in\ldo_{\devsp},\,\stm(\vf)=1}\frac{1}{\norm{\eps(\vf)}_{\yieldf}}\\
 & =\frac{1}{\inf_{\stm(\vf)=1}\norm{\eps(\vf)}_{\yieldf}}\end{align}
where we used $\sup_{x}(1/x)=1/\inf x$. We conclude that \begin{equation}
\frac{1}{\optimum_{\stm}}=\inf_{\stm(\vf)=1}\norm{\eps(\vf)}_{\yieldf}.\end{equation}
On the other hand, if we replace in the kinematic version of the limit
load the requirement $\st\in B$ by the requirement $\norm{\st}_{\yieldf}\les\yieldf_{0}$
which is the analog in our setting, we obtain \begin{align}
\lambda^{*} & =\inf_{\stm(\vf)=1}\braces{\sup_{\norm{\st}_{\yieldf}\les\yieldf_{0}}\braces{\eps^{*}(\st)(\vf)}}\\
 & =\inf_{\stm(\vf)=1}\braces{\sup_{\norm{\st}_{\yieldf}\les\yieldf_{0}}\braces{\st(\eps(\vf))}}.\end{align}
Using $\norm{\eps(\vf)}=\sup_{\norm{\st}\les1}\abs{\st(\eps(\vf))}$,
%
{} we have\begin{equation}
\lambda^{*}=\yieldf_{0}\inf_{\stm(\vf)=1}\norm{\eps(\vf)}_{\yieldf},\end{equation}
so indeed $\lambda^{*}=\yieldf_{0}/\optimum_{\stm}$.
\end{rem}

\bigskip{}

\noindent \textbf{\textit{Acknowledgments.}} This work was partially
supported by the Paul Ivanier Center for Robotics Research and Production
Management at Ben-Gurion University.

\appendix

\section{The Technical Lemma\label{sec:The-Technical-Lemma}}

\begin{lem}
Let $\vs$ and $\avs$ be two Banach spaces and $\fun\from\vs\tto\avs$
an isometric injection.

\noindent \textup{(}i\textup{)}~For each $\stm\in\vs^{*}$ there
is some (non-unique) $\st\in\avs^{*}$, such that\begin{equation}
\stm=\fun^{*}(\st).\end{equation}
\textup{(}ii\textup{)}~The dual norm of $\stm$ satisfies\begin{equation}
\norm{\stm}=\inf_{\st}\norm{\st},\end{equation}
where the infimum is taken over all $\st$ representing $\stm$, i.e.,
those satisfying\begin{equation}
\stm=\fun^{*}(\st).\end{equation}
\textup{(}iii\textup{)~}There is a $\widehat{\st}\in\avs^{*}$ such
that \begin{equation}
\norm{\stm}=\inf_{\st}\norm{\st}=\norm{\widehat{\st}}.\end{equation}

\end{lem}
\begin{proof}
Given $\stm\in\vs^{*}$, we may use the fact that \begin{equation}
\fun^{-1}\from\image\fun\subset\avs\tto\vs\end{equation}
is a well defined linear isometry to write\begin{equation}
\abs{\stm(\fun^{-1}(v))}\les\norm{\stm}\norm{\fun^{-1}(v)}=\norm{\stm}\norm{v}.\end{equation}
It follows that \begin{equation}
\stm\comp\fun^{-1}\from\image\fun\tto\reals\end{equation}
is a bounded linear functional on the subspace $\image\fun\subset\avs$.
We recall that the Hahn-Banach theorem asserts that if $\mathbf{U}\subset\avs$
is a vector subspace and $\tau$ is a bounded linear functional on
$\mathbf{U}$, then, $\tau$ may be extended to a bounded linear functional
$\st$ on $\avs$ such that\begin{equation}
\st(u)=\tau(u),\quad\text{for all}\quad u\in\mathbf{U},\end{equation}
and\begin{equation}
\norm{\st}=\sup_{v\in\mathbf{\avs}}\frac{\abs{\st(v)}}{\norm{v}}=\sup_{u\in\mathbf{\mathbf{U}}}\frac{\abs{\tau(u)}}{\norm{u}}=\norm{\tau}.\label{eq:HahnBanachNorm}\end{equation}
Applying the Hahn-Banach theorem to the situation at hand, we conclude
that the functional $\stm\comp\fun^{-1}$ may be extended to a linear
functional $\st$ on $\avs$ such that\begin{equation}
\st(u)=\stm\comp\fun^{-1}(u),\quad\text{for all}\quad u\in\image\fun,\end{equation}
or equivalently, \begin{equation}
\stm(\vf)=\st(\fun(\vf)).\end{equation}
By the definition of the dual mapping we conclude that $\stm=\fun^{*}(\st)$.

In general, for any $\st\in\avs^{*}$\begin{equation}
\norm{\fun^{*}(\st)}=\sup_{\vf\in\vs}\frac{\abs{\fun^{*}(\st)(\vf)}}{\norm{\vf}}=\sup_{\vf\in\vs}\frac{\abs{\st(\fun(\vf))}}{\norm{\vf}}\les\sup_{\vf\in\vs}\frac{\norm{\st}\norm{\fun(\vf)}}{\norm{\fun(\vf)}}\,,\end{equation}
so\begin{equation}
\norm{\fun^{*}(\st)}\les\norm{\st}.\end{equation}
On the other hand, for any $\st\in\avs^{*}$, such that $\stm=\fun^{*}(\st)$\begin{equation}
\sup_{\vf\in\vs}\frac{\abs{\stm(\vf)}}{\norm{\vf}}=\sup_{v\in\image\fun}\frac{\abs{\stm\comp\fun^{-1}(v)}}{\norm{v}}=\sup_{v\in\image\fun}\frac{\abs{\st(v)}}{\norm{v}},\end{equation}
so by the Hahn-Banach theorem\begin{equation}
\norm{\stm}=\norm{\fun^{*}(\st)}=\norm{\widehat{\st}},\end{equation}
where $\widehat{\st}$ is the element of $\avs^{*}$ extending $\stm\comp\fun^{-1}$
and having the same norm as in Equation (\ref{eq:HahnBanachNorm}).

We conclude that \begin{equation}
\norm{\stm}=\inf_{\st}\norm{\st},\end{equation}
where infimum is taken over all $\st$ satisfying $\stm=\fun^{*}(\st)$.
The infimum is attained for $\widehat{\st}$ as above.
\end{proof}


\begin{thebibliography}{10}
\bibitem{Chris1}Christiansen, E., Limit analysis in plasticity as
a mathematical programming problem, \emph{Calcolo}, \textbf{17}, 41--65
(1980).

\bibitem{Chris2}Christiansen, E., On the collapse solution in limit
analysis, \emph{Archive for Rational Mechanics and Analysis,} \textbf{91},
119--135 (1986).

\bibitem{PeretzSegev05}Peretz, R., Segev, R., Bounds on the trace
mapping of $\ld$-fields, accepted for publication, \emph{Computers
and Mathematics with Applications,} arXiv:Math.AP/0505006 (2005).

\bibitem{segev03}Segev, R., Generalized stress concentration factors,
\emph{Mathematics and Mechanics of Solids,} doi: 10.1177/1081286505044131
(2005).

\bibitem{segev04}Segev, R., Generalized stress concentration factors
for equilibrated forces and stresses, accepted for publication, \emph{Journal
of Elasticity}, \textbf{81}, 293 -- 315 (2005).

\bibitem{Taylor}Taylor, A.E., \emph{Introduction to Functional Analysis},
Wiley, (1958).

\bibitem{TemamStrang1}Temam, R., Strang, G., Functions of bounded
deformations, \emph{Archive for Rational Mechanics and Analysis},
\textbf{75}, 7--21 (1980).

\bibitem{TemamStrang2}Temam, R., Strang, G., Duality and relaxation
in the variational problems of plasticity, \emph{J de Mecanique,}
\textbf{19}, 493--527 (1980).

\bibitem{Temam81}Temam, R., On the continuity of the trace of vector
functions with bounded deformation, \textit{Applicable Analysis},
\textbf{11}, 291--302 (1981).

\bibitem{Temam85}Temam, R., \emph{Mathematical Problems in Plasticity,}
(a translation of \emph{Problemes mathematiques en plasticite}, Bordas,
Paris, 1983) Gauthier-Villars, Paris (1985).
\end{thebibliography}
\end{document}